\documentclass[10pt]{article}

\usepackage{algorithm2e}
\usepackage{xcolor}

\SetAlFnt{\footnotesize}

\SetAlCapSty{xAlCapSty}

\SetCommentSty{xCommentSty}

\SetNlSty{mynlfont}{}{} 

\LinesNumbered

\SetSideCommentRight

\DontPrintSemicolon

\RestyleAlgo{algoruled}

\usepackage[thinlines]{easytable}

\usepackage{color}      % use if color is used in text
\usepackage{placeins}
\usepackage{indentfirst}
\usepackage{multirow}
\usepackage{amssymb}
\usepackage{amsmath,latexsym}
\usepackage{url}
\usepackage{enumerate}
\usepackage{graphicx}

\usepackage{booktabs}

\usepackage{subcaption}
\usepackage{cleveref}
\usepackage{mathtools}

\begin{document}

\title{\bf Grid Quality Measures for Iterative Convergence}

 \author{
{Hiroaki Nishikawa}\thanks{Associate Research Fellow
({hiro@nianet.org}),
%{National Institute of Aerospace},s
 100 Exploration Way, Hampton, VA 23666 USA,  }\\
  {\normalsize\itshape 
{National Institute of Aerospace},
 Hampton, VA 23666, USA}
}

%\author{Hiroaki Nishikawa and Yoshitaka Nakashima}
%\date{February 28, 2019}
\date{\today}
\maketitle
%%%%%%%%%%%%

\begin{abstract} 
In this paper, we discuss two grid-quality measures, F- and G-measures, in relation to iterative convergence of an implicit unstructured-grid Navier-Stokes solver.
The F-measure is a lower bound of a least-squares gradient, which is a purely geometrical quantity defined in each cell and thus can be computed for a given grid: 
faster convergence is expected for a grid with a lower value of the F-measure. The G-measure is a least-squares gradient of a specified function 
around each cell, with the minimum value of zero. Faster convergence is expected for a smaller value of the G-measure towards zero. 
In this paper, we investigate these measures for inviscid and viscous problems with unstructured grids in two dimensions. 
\end{abstract}
 
%%%%%%%%%%%%%%%%%%%%%%%%%%%%%%%%%%%%%%%%%%%%%%%%%%%%%%%%%%%
\section{Introduction}
%%%%%%%%%%%%%%%%%%%%%%%%%%%%%%%%%%%%%%%%%%%%%%%%%%%%%%%%%%
 
The quality of an unstructured grid has a significant impact on accuracy and robustness of practical computational fluid
dynamics solvers \cite{VeluriRoyLuke:CF2012,KallinderisFotia:JCP2015,YouMittalWangMoin:JCP2006,DennerWachem:JCP2015,ghoreyshi_etal:AIAA2015-0407,nc_cc_comparison_viscous:AIAAJ2010,Dannenhoffer:AIAA2012-0610,pointwise_article_link}. Various measures are available to assess the grid quality: grid skewness, aspect ratio, stretching factor, face angles, etc., and each of them has some implications to accuracy and robustness of an unstructured-grid solver. For example, explicit schemes tend to become less accurate and can get unstable on skewed structured grids \cite{YouMittalWangMoin:JCP2006} while implicit solvers tend to slow down on highly skewed unstructured grids \cite{diskin_thomas:AIAAJ2011,nishikawa_nakashima_watanabe:jcp2017}. These grid-quality measures can be used to assess the quality of general unstructured grids, but they are usually too general to predict the performance of a particular discretization. In order to gain accurate indications for a given solver, we need measures closely related to the target solver. One such example is the grid skewness measure defined at a face shared by two adjacent cells, say, 1 and 2, as the dot product of the unit face normal vector pointing from 1 to 2 and the unit vector pointing from the centroid of the cell 1 to that of the cell 2. As discussed in Ref.\cite{nishikawa_stencil:JCP2019}, this particular skewness measure explicitly appears in the finite-volume discsretization of the viscous terms and thus has significant impact on the high-frequency damping property and implicit solver convergence. It was shown to have some impact on the inviscid schemes as well but it is difficult to show clearly how the grid skewness measure affects inviscid discretizations.   
In this paper, inspired by the novel least-squares (LSQ) gradient stencil techniques proposed in Ref.\cite{nishikawa_stencil:JCP2019}, we introduce grid-quality measures in relation to iterative convergence of an implicit Navier-Stokes solver, focusing on a finite-volume inviscid discretization. 

In Ref.\cite{nishikawa_stencil:JCP2019}, novel gradient-stencil construction techniques are proposed for stable iterative convergence of an implicit finite-volume solver based on two observations: a solver is expected to be more stable (1) when a gradient stencil is symmetric (e.g., structured grids) and (2) when solution gradients are smaller in magnitude (e.g., a first-order scheme). In the latter, a lower bound of a LSQ gradient can be derived, which can be computed for a given grid. Ref.\cite{nishikawa_stencil:JCP2019} suggests that it can be used to select extra neighbors and adde them to a gradient stencil that will reduce the lower bound. Numerical results show that the resulting gradient stencils successfully lead to more stable convergence. In this paper, we investigate the use of the lower bound as a grid-quality measure for robust iterative convergence. That is, we investigate whether the lower bound computed for a given grid can indicate the performance of the implicit finite-volume solver on that grid. The lower bound, however, is not bounded and thus the minimum value is not known. It is still useful as a smaller value indicates better iterative convergence, but it would be reasonable to consider an alternative measure that has the minimum (i.e., the best) value. For this purpose, we introduce the magnitude of the gradient of a known function, which should be zero on an isotropic or symmetric stencil. 

\begin{figure}[htbp!]
    \centering
          \begin{subfigure}[t]{0.23\textwidth}
        \includegraphics[width=\textwidth,trim=3 3 3 3,clip]{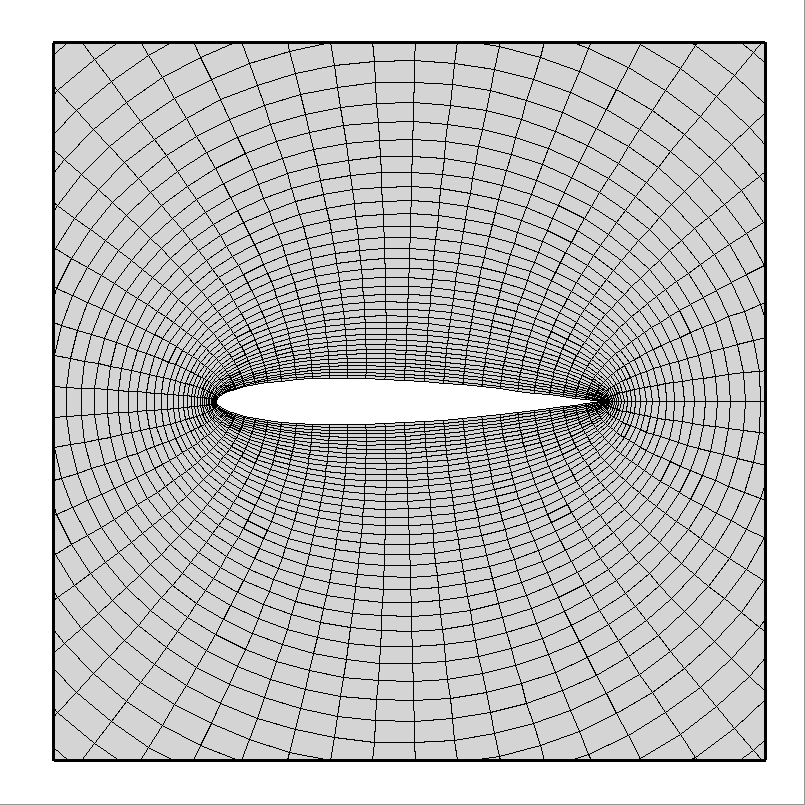}
                  \caption{Structured.}
          \label{fig:airfoil_str}
      \end{subfigure}
    \hfill    
          \begin{subfigure}[t]{0.23\textwidth}
        \includegraphics[width=\textwidth,trim=3 3 3 3,clip]{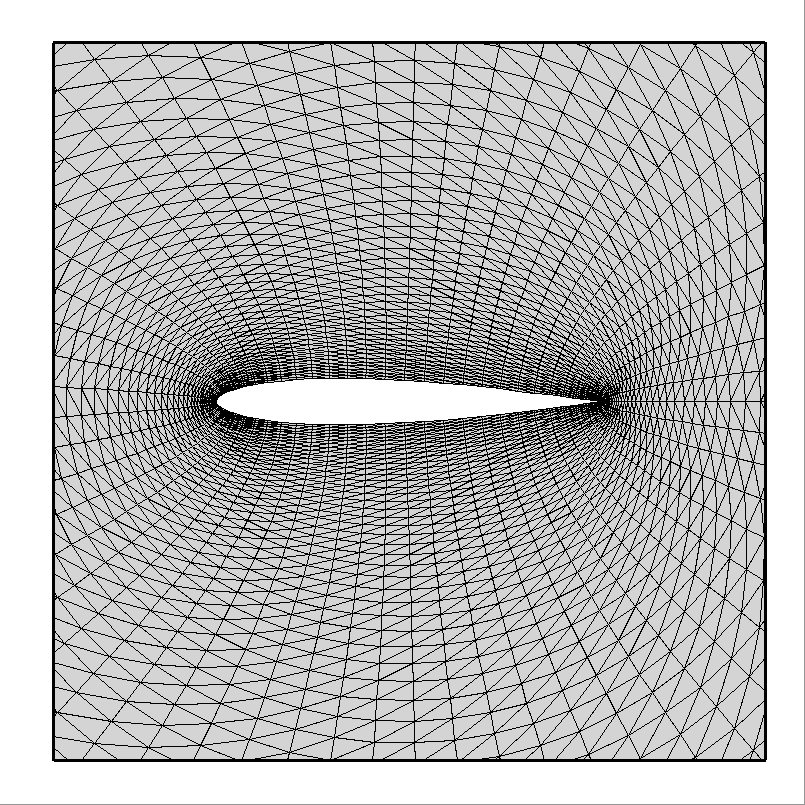}
                  \caption{Triangular.}
          \label{fig:airfoil_tri}
      \end{subfigure}
    \hfill    
                \begin{subfigure}[t]{0.23\textwidth}
        \includegraphics[width=\textwidth,trim=3 3 3 3,clip]{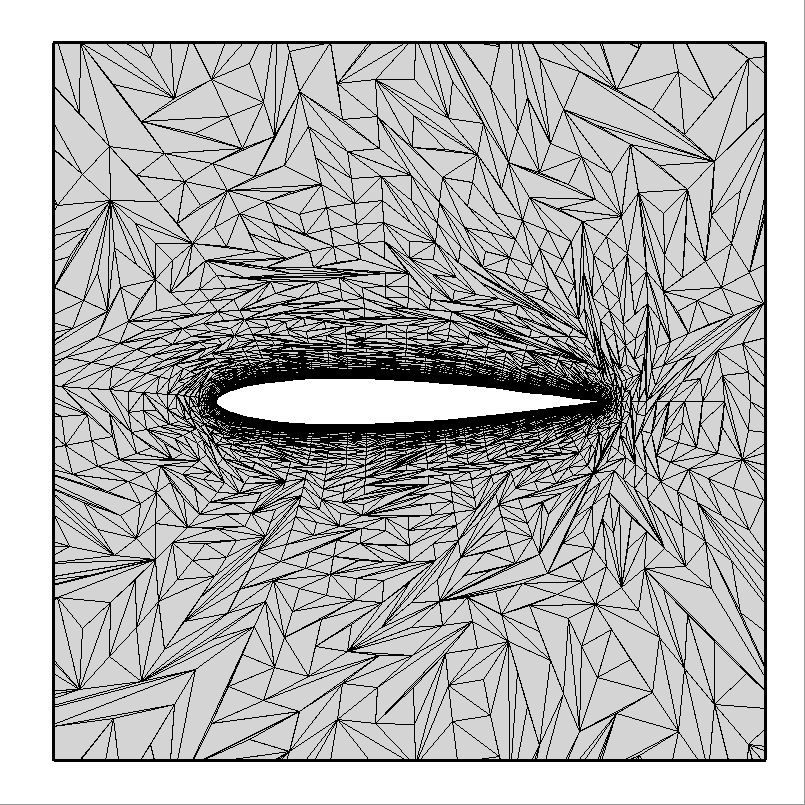}
                  \caption{Irregular.}
          \label{fig:airfoil_irrg}
      \end{subfigure}
            \caption{
\label{fig:airfoil_grids_introduction}%
Different types of grids.
} 
\end{figure}

Note that one could make fairly good prediction if a grid can be inspected beforehand. For example, see grids shown in Figure \ref{fig:airfoil_grids_introduction}.
Clearly, we would expect that a finite-volume solver would converge faster on the structured grid in Figure \ref{fig:airfoil_str} and have difficulties on the irregular grid 
\ref{fig:airfoil_irrg}. However, such requires the visualization of grid in each practical application and the knowledge and experience of a user. It is highly desirable to remove such a process, especially in automated simulations with grid adaptation \cite{Kleb_etal_aiaa2019-2948}. In this paper, we will focus on iterative convergence of an implicit finite-volume solver and investigate the use of the proposed measures as an indication of the iterative convergence. A cell-centered finite-volume method will be considered, but the proposed measures can be easily applied to node-centered methods.

The paper is organized as follows.
In Section 2, we describe the lower bound, called the F-measure, and the magnitude of the gradient of a known function, called the G-measure. 
In Section 3, we present numerical results to investigate how these measures predict the iterative convergence of an implicit finite-volume solver.
In Section 4, we conclude the paper with remarks.

%%%%%%%%%%%%%%%%%%%%%%%%%%%%%%%%%%%%%%%%
%%%%%%%%%%%%%%%%%%%%%%%%%%%%%%%%%%%%%%%%%

%%%%%%%%%%%%%%%%%%%%%%%%%%%%%%%%%%%%%%%%%%%%%%%%%%%%%%%%%%%
\section{Grid Quality Measures}
%%%%%%%%%%%%%%%%%%%%%%%%%%%%%%%%%%%%%%%%%%%%%%%%%%%%%%%%%%

We consider an implicit defect-correction solver, which uses the exact Jacobian of the low-order residual (i.e., zero LSQ gradients)
for solving a system of second-order cell-centered finite-volume residual equations for the Navier-Stokes equations. See Ref.\cite{nishikawa_stencil:JCP2019} for
details of the solver and the discretization. The implicit solver becomes Newton's method, which is fast and robust for solving a system of low-order residual equations. It implies that the solver converges faster and is more robust if the LSQ gradients are closer to zero. Based on this observation, a gradient stencil
augmentation technique was devised \cite{nishikawa_stencil:JCP2019}, where extra cells are added to a gradient stencil that reduce the lower bound of the LSQ gradient. 
The lower bound is given by the F-measure as described below.

\subsection{F-Measure}
Consider the normal equation for a linear LSQ gradient method for computing the gradient $( \overline{\partial_x u}_j , \overline{\partial_y u}_j )$ for a variable $u$:
\begin{eqnarray}
{\bf A}^T {\bf A} {\bf x} = {\bf A}^T {\bf b}, 
\label{lsq_normal_system}
\end{eqnarray}
where
\begin{eqnarray}
{\bf A}^T {\bf A} = 
\left[
\begin{array}{cc}
\displaystyle  \sum_{k \in \{ g_j \}} w_k^2 \Delta x_{k}^2 & \displaystyle  \sum_{k \in \{ g_j \}} w_k^2 \Delta x_{k} \Delta y_{k} \\ [4.5ex]
\displaystyle   \sum_{k \in \{ g_j \}} w_k^2 \Delta x_{k} \Delta y_{k}  &\displaystyle    \sum_{k \in \{ g_j \}} w_k^2 \Delta y_{k}^2
\end{array}
\right], 
\quad
{\bf x} = 
\left[
\begin{array}{c}
\displaystyle  \overline{\partial_x u}_j \\ [2.5ex]
\displaystyle  \overline{\partial_y u}_j
\end{array}
\right], 
\quad
 {\bf A}^T {\bf b} = 
\left[
\begin{array}{c}
 \displaystyle   \sum_{k \in \{ g_j \}} w_k^2 \Delta x_{k}  \Delta u_{k} \\ [4.5ex]
\displaystyle   \sum_{k \in \{ g_j \}} w_k^2 \Delta y_{k}  \Delta u_{k}
\end{array}
\right],
\end{eqnarray} 
\begin{eqnarray}
 \Delta x_{k} = x_k - x_j , \quad \Delta y_{k} = y_k - y_j , \quad  \Delta u_{k} = u_k - u_j, \quad k = 1, 2, , \cdots, N,
\end{eqnarray}
where $ \{ g_j \}$ denotes a set of neighbor cells of the cell $j$, $N$ is the number of the neighbors, $w_k$'s are weights, $(x_j,y_j)$ is the centroid coordinates of the cell $j$, and
similarly for the neighbor cells. 
As derived in Ref.\cite{nishikawa_stencil:JCP2019}, the lower bound of the gradient is proportional to the following quantity:
\begin{eqnarray}
 F = \frac{  s  }{ || {\bf A}^T {\bf A} ||_F },  \quad
s = \sum_{k \in \{ g_j \}}   w_k^2  d_k, \quad
w_k = \frac{1}{ d_k^p}, \quad
d_k = \sqrt{ (x_k - x_j )^2 +  (y_k - y_j )^2  },
\end{eqnarray}
where $p = 0$ for the inviscid reconstruction and $p=1$ for the viscous reconstruction \cite{nishikawa_stencil:JCP2019,WhiteNishikawaBaurle_aiaa2019-0127}.
Although it is merely the lower bound, it has been found useful for constructing a LSQ stencil for 
stable iterative convergence on distorted grids: i.e., more stable with gradient stencils of lower value of $F$ \cite{nishikawa_stencil:JCP2019}. 
The measure $F$ is referred here to as the F-measure.
Note that the value of $p$ should be chosen in consistent with the value used in the actual residual calculation. In this study, we focus on the inviscid discretization,
and therefore set $p=0$. Another point to note about the F-measure is that its minimum value is not clearly known and thus it is difficult to assess a given grid without reference to other data. 

\subsection{G-Measure} 
To estimate the gradient magnitude more directly and with a known minimum value of a measure, we propose to compute 
the LSQ gradient of a known function, e.g., a Gaussian function, and use its magnitude as a grid-quality measure.
To eliminate ambiguity as much as possible, we employ an isotropic function in normalized coordinates: 
\begin{eqnarray}
u = \exp( -(\tilde{x}^2+\tilde{y}^2) ), \quad
\tilde{x} = (x - x_j) / s_{max}, \quad \tilde{y} =  (y - y_j) / s_{max},
\end{eqnarray}
where $s_{max} = \sqrt{ x_{max}^2 + y_{max}^2 } $, $(x_{max}, y_{max})$ is the centroid coordinates of the neighbor cell farthest from the target cell $j$, 
and compute the LSQ gradient:
\begin{eqnarray}
( \, \overline{\partial_x u}_j , \overline{\partial_y u}_j  \, )
= \sum_{k \in \{ g_j \} } ( c^x_k, c^y_k ) \left[ \, 
 \exp( -(\tilde{x}_k^2+\tilde{y}_k^2) ) -  \exp( -(\tilde{x}_j^2+\tilde{y}_j^2) )
\, \right],
\label{lsq_grad}
\end{eqnarray}
where $ (  c^x_k, c^y_k) $ are the LSQ gradient coefficients \cite{nishikawa_stencil:JCP2019}.
As a grid-quality measure, we define 
\begin{eqnarray}
G = \sqrt{  (\overline{\partial_x u}_j)^2 +  (\overline{\partial_y u}_j)^2 }. 
\label{lsq_grad}
\end{eqnarray}
The exact gradient is $(0,0)$ at the target cell centroid, and therefore the ideal value of $G$ is zero. 
Our expectation is that the implicit solver is more robust and converges faster if this measure is closer to zero.
We call this measure the G-measure. 

\subsection{Remarks} 

In this study, we compute the minimum, maximum, and average of these measures over all cells, and 
investigate how they indicate the performance of iterative convergence of an implicit solver. 
Local values of the F- and G-measures may be used for grid adaptation (e.g., swap an edge if it will reduce the G-measures 
over the cells affected by the swapping). Such is left as future work. 

%%%%%%%%%%%%%%%%%%%%%%%%%%%%%%%%%%%%%%%%%%%%%%%%%%%%%%%%%%%
%%%%%%%%%%%%%%%%%%%%%%%%%%%%%%%%%%%%%%%%%%%%%%%%%%%%%%%%%%%
\section{Numerical Results}
\label{new_stencil}
%%%%%%%%%%%%%%%%%%%%%%%%%%%%%%%%%%%%%%%%%%%%%%%%%%%%%%%%%%%
%%%%%%%%%%%%%%%%%%%%%%%%%%%%%%%%%%%%%%%%%%%%%%%%%%%%%%%%%%%

In all numerical results, the implicit defect-correction solver as described in Ref.\cite{nishikawa_stencil:JCP2019} is employed 
for solving the Euler and Navier-Stokes equations in two dimensions. The objective is to investigate whether the implicit solver will converge faster for grids with smaller F- or G-measures as expected.

%%%%%%%%%%%%%%%%%%%%%%%%%%%%%%%%%
\subsection{Viscous manufactured solution at $M_\infty = 0.2$ and $Re_\infty = 40$}
%%%%%%%%%%%%%%%%%%%%%%%%%%%%%%%%%
%%%%%%%%%%%%%%%%%%%%%%%%%%%%%%%%

\begin{figure}[htbp!]
    \centering
          \begin{subfigure}[t]{0.39\textwidth}
        \includegraphics[width=\textwidth,trim=2 2 2 2,clip]{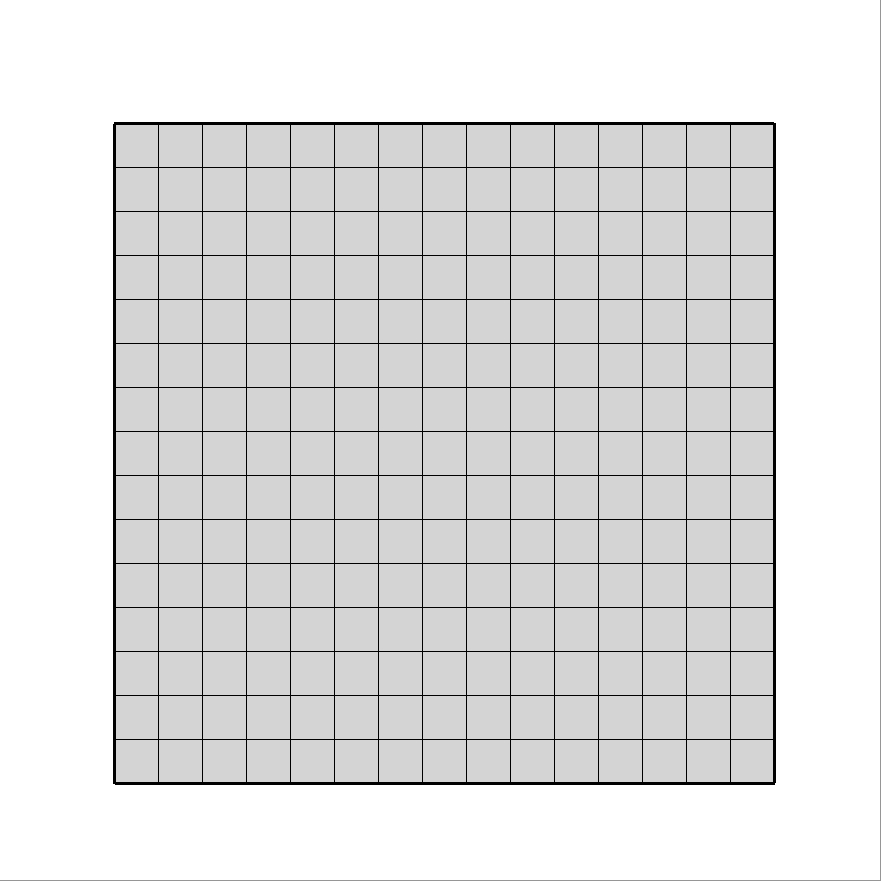}
                  \caption{Grid00.}
          \label{fig:mms_grid00}
      \end{subfigure}
    \hfill    
          \begin{subfigure}[t]{0.39\textwidth}
        \includegraphics[width=\textwidth,trim=2 2 2 2,clip]{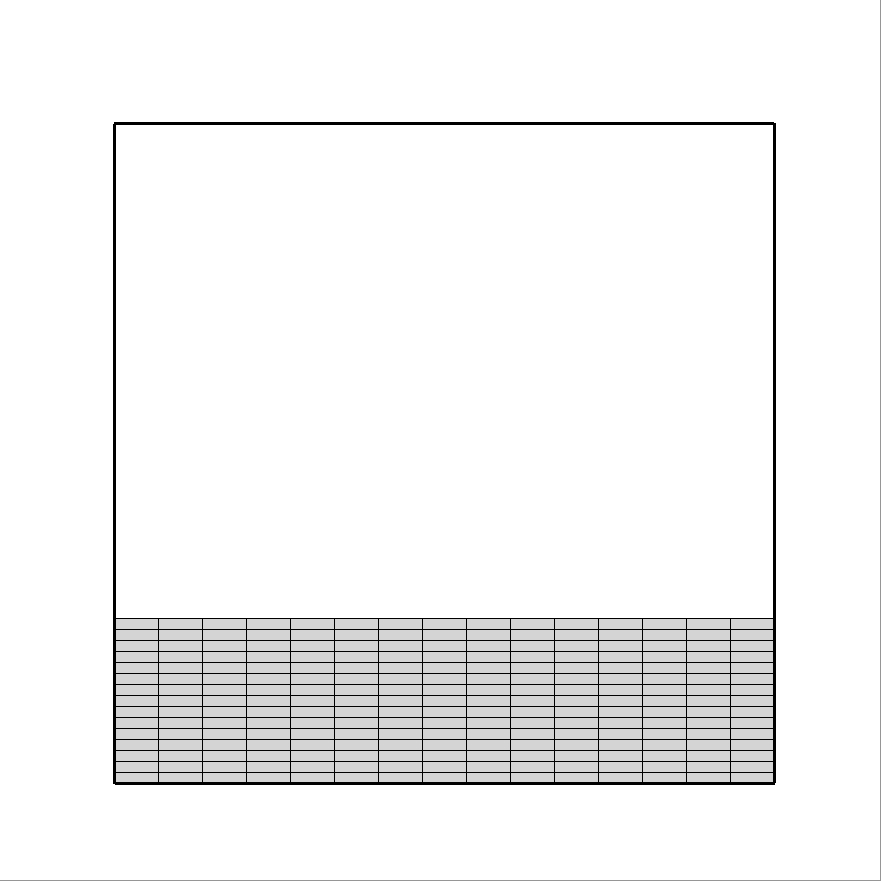}
                  \caption{Grid01.}
          \label{fig:mms_grid01}
      \end{subfigure}
    \hfill    
          \begin{subfigure}[t]{0.39\textwidth}
        \includegraphics[width=\textwidth,trim=2 2 2 2,clip]{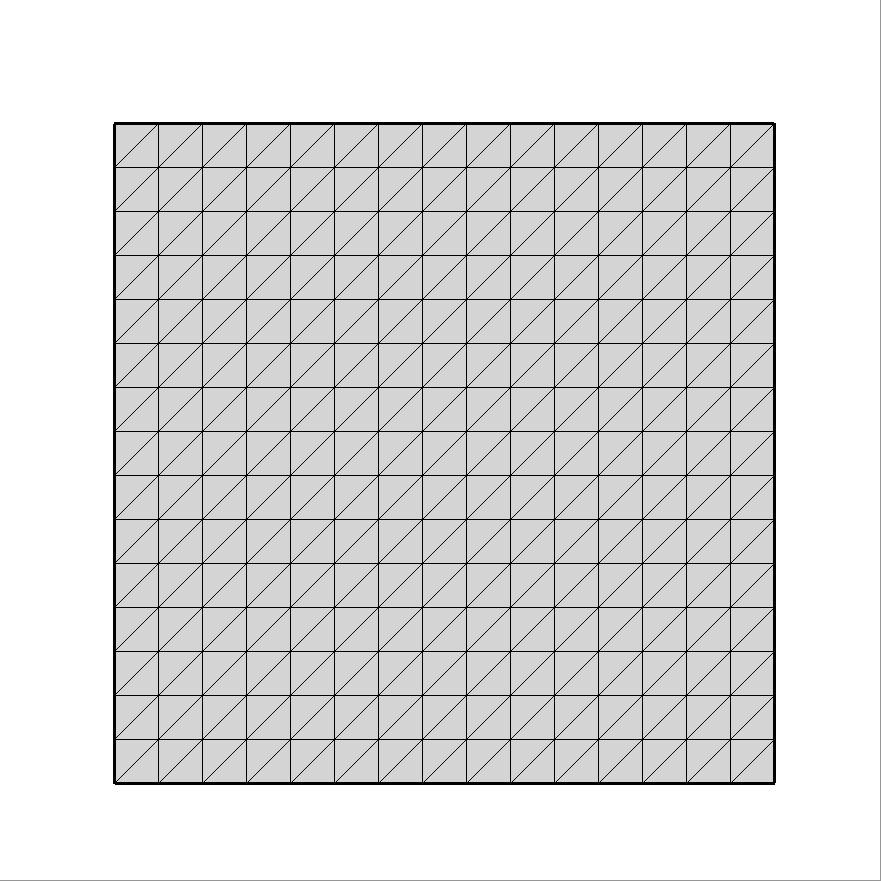}
                  \caption{Grid02.}
          \label{fig:mms_grid02}
      \end{subfigure}
    \hfill    
          \begin{subfigure}[t]{0.39\textwidth}
        \includegraphics[width=\textwidth,trim=2 2 2 2,clip]{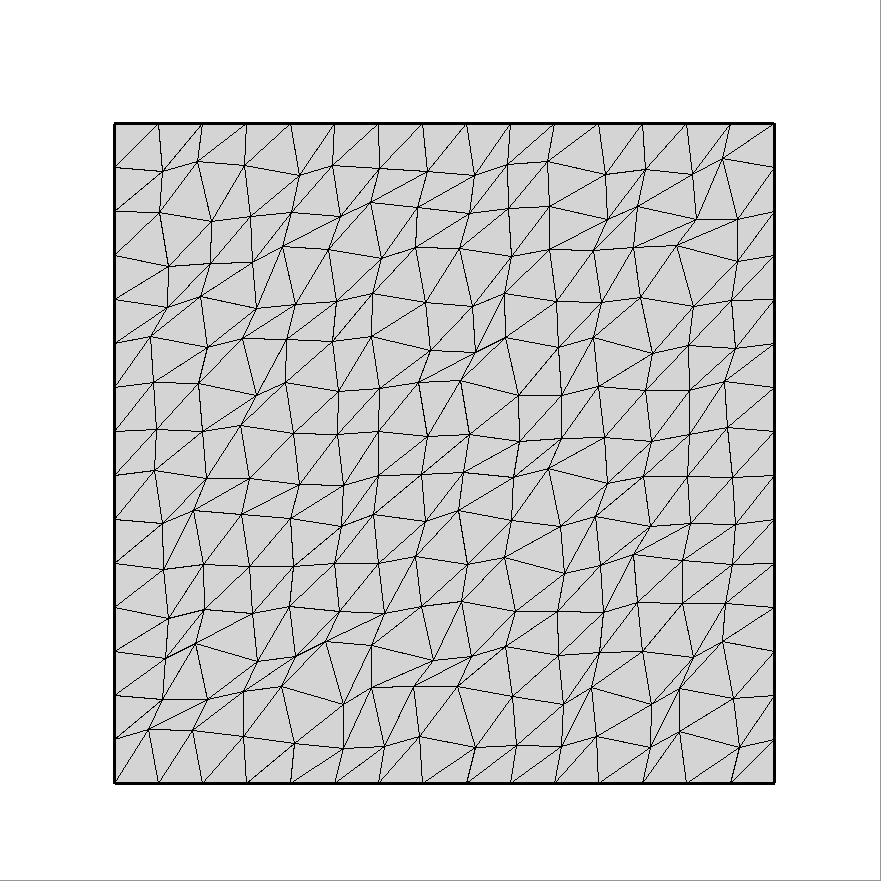}
                  \caption{Grid10.}
          \label{fig:mms_grid10}
      \end{subfigure}
            \caption{
\label{fig:mms_grids}%
Coarsest grids (level 1) for the MMS test case.
} 
\end{figure}

\begin{figure}[htbp!]
    \centering
          \begin{subfigure}[t]{0.4\textwidth}
        \includegraphics[width=\textwidth,trim=0 0 0 0,clip]{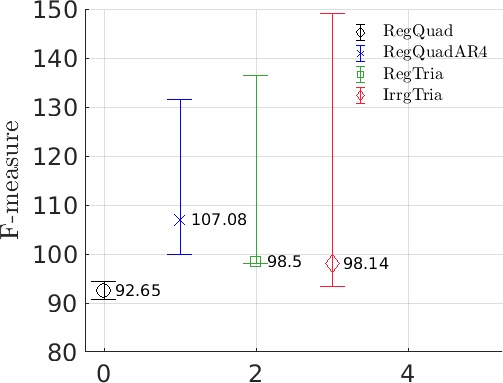}
          \caption{F-measure.}
          \label{fig:mms_f}
      \end{subfigure}
          \hfill    
          \begin{subfigure}[t]{0.4\textwidth}
        \includegraphics[width=\textwidth,trim=0 0 0 0 ,clip]{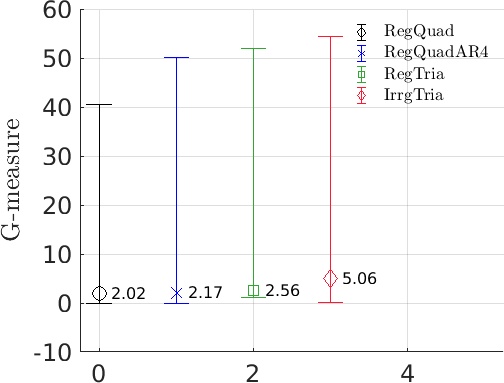}
                  \caption{G-measure.}
          \label{fig:mms_g}
      \end{subfigure}
                  \caption{
\label{fig:mms_f_and g}%
F- and G-measures for the finest (level 5) grids.
} 
      \hfill    
\end{figure}

We consider solving the Navier-Stokes equations for manufactured solutions \cite{nishikawa_centroid:JCP2020} with $T_\infty = 300 [K]$
and the Prandtl number 0.72 on four different types of grids as shown in Figure \ref{fig:mms_grids}: quadrilateral grid, quadrilateral grid with aspect ratio of 4, regular triangular grid, and irregular triangular gridd. Each grid has five levels: 16$\times$16, 32$\times$32, 48$\times$48, 64$\times$64, 80$\times$80 nodes. 
As mentioned earlier, it may be easy to predict the performance of the solver by inspecting these grids, but our interest is to investigate whether the F- and G-measures can serve as accurate quantitative measure for the solver performance. The F- and G-measures computed on the finest grid are shown in Figure \ref{fig:mms_f_and g}. It is noted that the G-measure is more consistent with our intuition that the solver converges well in the order for Grid00, Grid01, Grid 03, and then Grid10. 

\begin{figure}[htbp!]
    \centering  
          \begin{subfigure}[t]{0.24\textwidth}
        \includegraphics[width=\textwidth,trim=2 2 2 2,clip]{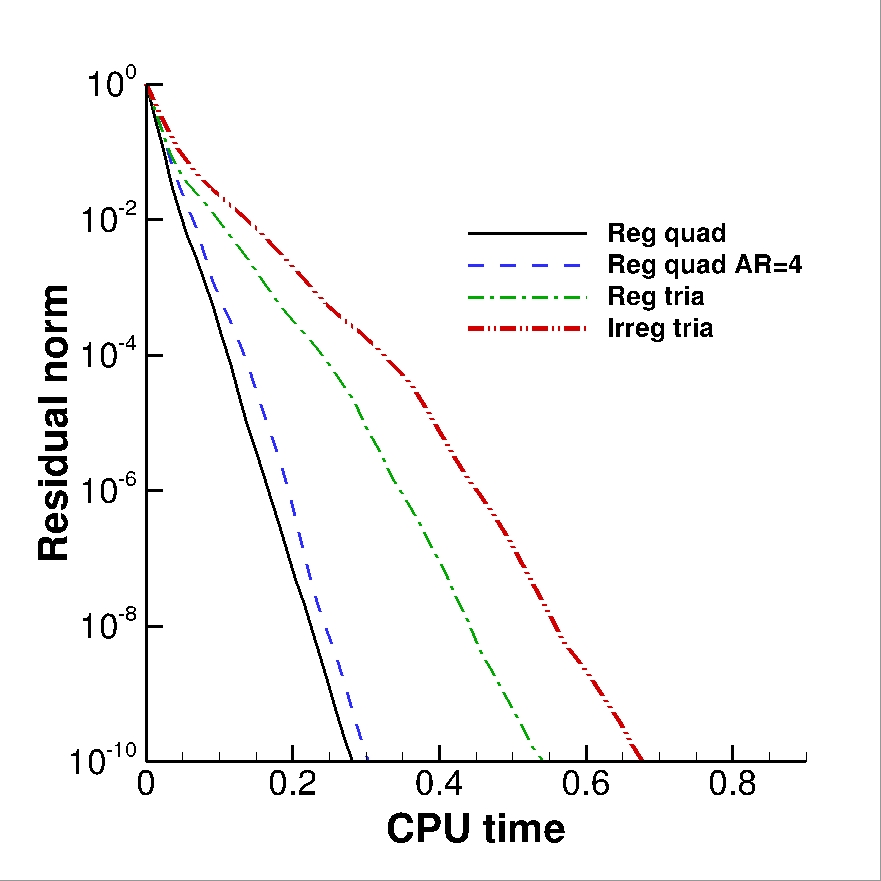}
          \caption{Grid00.}
      \end{subfigure}
          \hfill    
          \begin{subfigure}[t]{0.24\textwidth}
        \includegraphics[width=\textwidth,trim=2 2 2 2,clip]{./F_measure_visc_mms/figures/figs_res1_vs_cpu_g01.jpeg}
          \caption{Grid01.}
      \end{subfigure}
          \begin{subfigure}[t]{0.24\textwidth}
        \includegraphics[width=\textwidth,trim=2 2 2 2,clip]{./F_measure_visc_mms/figures/figs_res1_vs_cpu_g01.jpeg}
          \caption{Grid02.}
      \end{subfigure}
          \hfill    
          \begin{subfigure}[t]{0.24\textwidth}
        \includegraphics[width=\textwidth,trim=2 2 2 2,clip]{./F_measure_visc_mms/figures/figs_res1_vs_cpu_g01.jpeg}
          \caption{Grid10.}
      \end{subfigure}
            \caption{
\label{fig:mms_res_cpu_glvl01}%
Grid level 1: Residual norm versus CPU time for the MMS test case.
} 
\end{figure}

\begin{figure}[htbp!]
    \centering  
          \begin{subfigure}[t]{0.24\textwidth}
        \includegraphics[width=\textwidth,trim=2 2 2 2,clip]{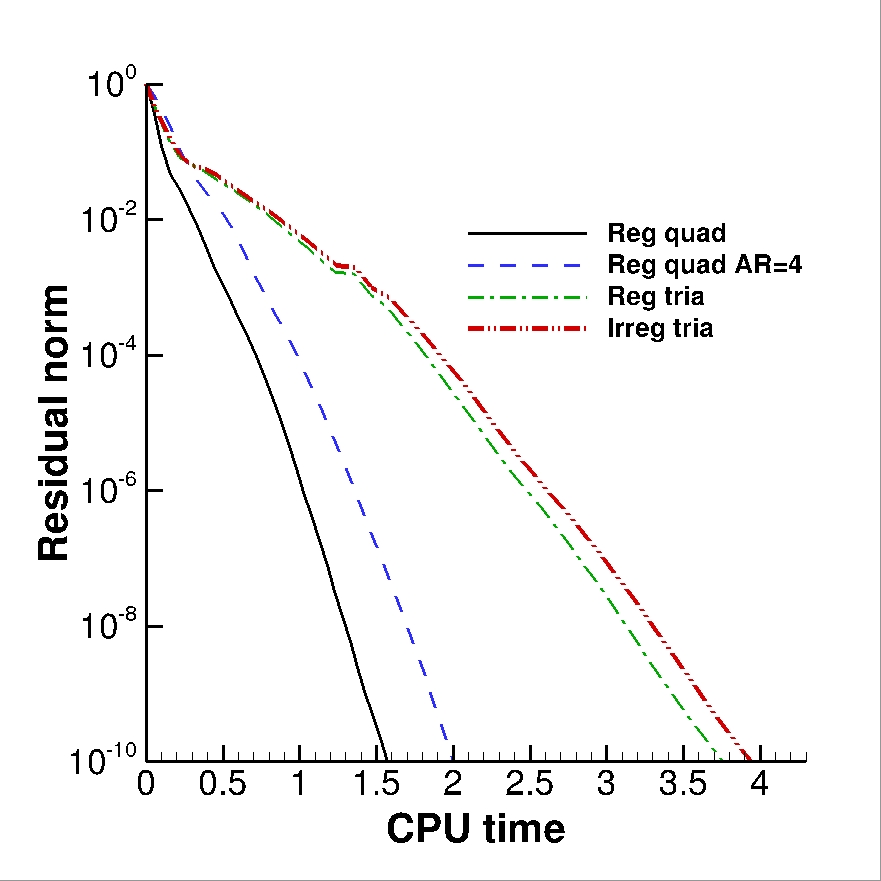}
          \caption{Grid00.}
      \end{subfigure}
          \hfill    
          \begin{subfigure}[t]{0.24\textwidth}
        \includegraphics[width=\textwidth,trim=2 2 2 2,clip]{./F_measure_visc_mms/figures/figs_res1_vs_cpu_g02.jpeg}
          \caption{Grid01.}
      \end{subfigure}
          \begin{subfigure}[t]{0.24\textwidth}
        \includegraphics[width=\textwidth,trim=2 2 2 2,clip]{./F_measure_visc_mms/figures/figs_res1_vs_cpu_g02.jpeg}
          \caption{Grid02.}
      \end{subfigure}
          \hfill    
          \begin{subfigure}[t]{0.24\textwidth}
        \includegraphics[width=\textwidth,trim=2 2 2 2,clip]{./F_measure_visc_mms/figures/figs_res1_vs_cpu_g02.jpeg}
          \caption{Grid10.}
      \end{subfigure}
            \caption{
\label{fig:mms_res_cpu_glvl02}%
Grid level 2: Residual norm versus CPU time for the MMS test case.
} 
\end{figure}

\begin{figure}[htbp!]
    \centering  
          \begin{subfigure}[t]{0.24\textwidth}
        \includegraphics[width=\textwidth,trim=2 2 2 2,clip]{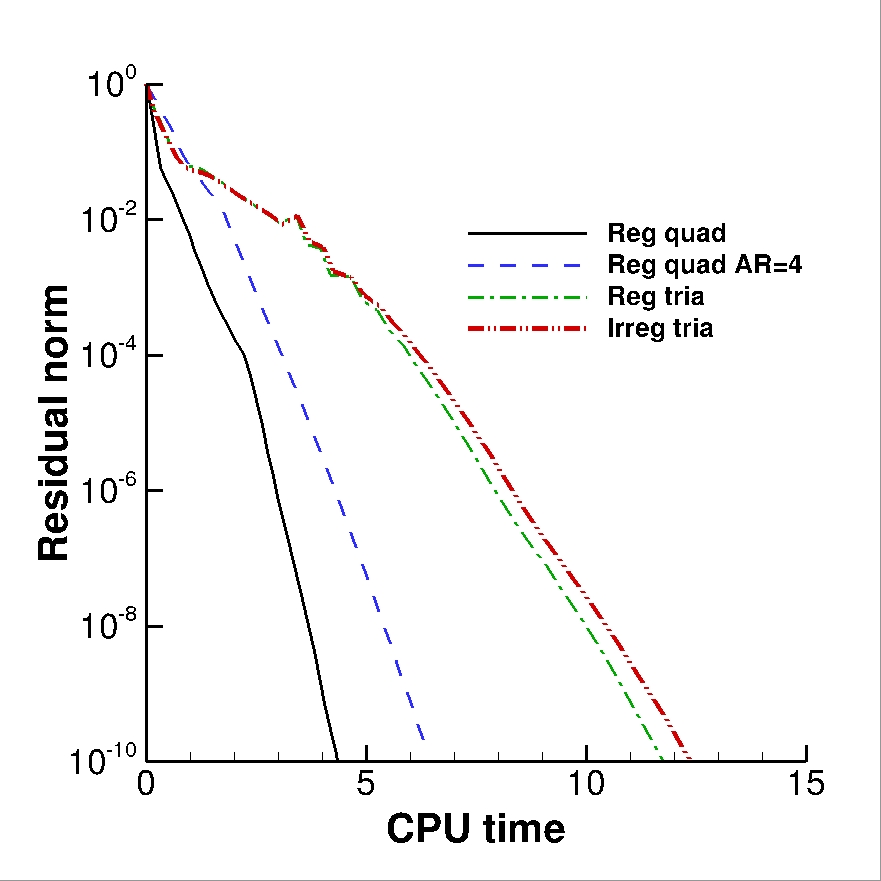}
          \caption{Grid00.}
      \end{subfigure}
          \hfill    
          \begin{subfigure}[t]{0.24\textwidth}
        \includegraphics[width=\textwidth,trim=2 2 2 2,clip]{./F_measure_visc_mms/figures/figs_res1_vs_cpu_g03.jpeg}
          \caption{Grid01.}
      \end{subfigure}
          \begin{subfigure}[t]{0.24\textwidth}
        \includegraphics[width=\textwidth,trim=2 2 2 2,clip]{./F_measure_visc_mms/figures/figs_res1_vs_cpu_g03.jpeg}
          \caption{Grid02.}
      \end{subfigure}
          \hfill    
          \begin{subfigure}[t]{0.24\textwidth}
        \includegraphics[width=\textwidth,trim=2 2 2 2,clip]{./F_measure_visc_mms/figures/figs_res1_vs_cpu_g03.jpeg}
          \caption{Grid10.}
      \end{subfigure}
            \caption{
\label{fig:mms_res_cpu_glvl03}%
Grid level 3: Residual norm versus CPU time for the MMS test case.
} 
\end{figure}

\begin{figure}[htbp!]
    \centering  
          \begin{subfigure}[t]{0.24\textwidth}
        \includegraphics[width=\textwidth,trim=2 2 2 2,clip]{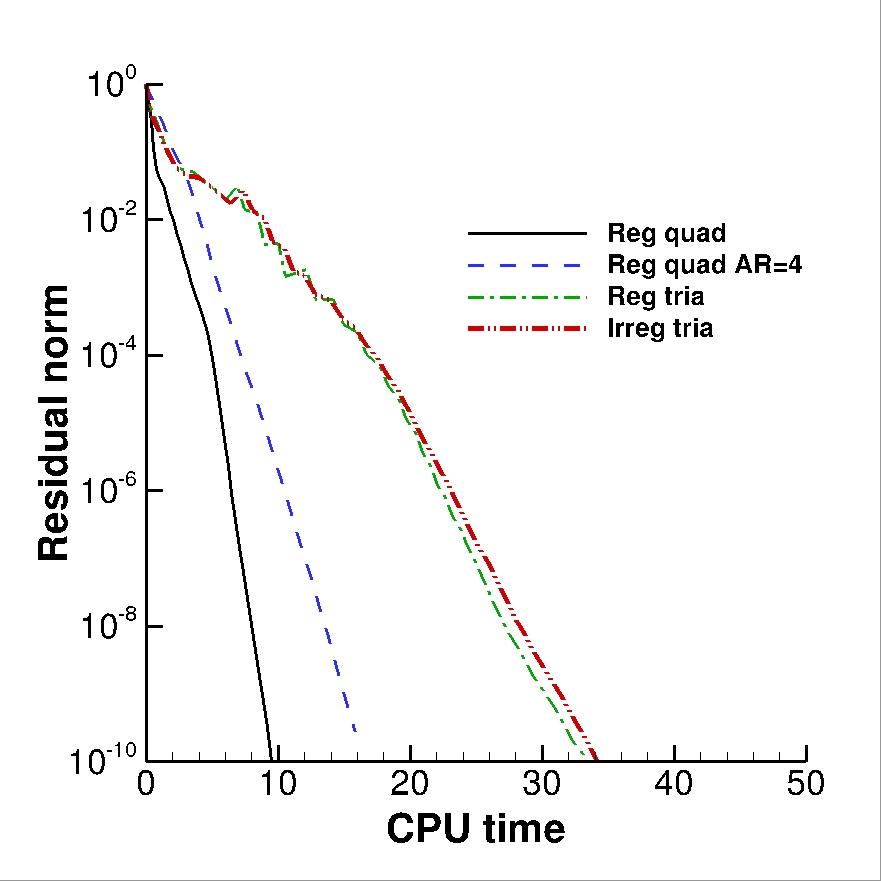}
          \caption{Grid00.}
      \end{subfigure}
          \hfill    
          \begin{subfigure}[t]{0.24\textwidth}
        \includegraphics[width=\textwidth,trim=2 2 2 2,clip]{./F_measure_visc_mms/figures/figs_res1_vs_cpu_g04.jpeg}
          \caption{Grid01.}
      \end{subfigure}
          \begin{subfigure}[t]{0.24\textwidth}
        \includegraphics[width=\textwidth,trim=2 2 2 2,clip]{./F_measure_visc_mms/figures/figs_res1_vs_cpu_g04.jpeg}
          \caption{Grid02.}
      \end{subfigure}
          \hfill    
          \begin{subfigure}[t]{0.24\textwidth}
        \includegraphics[width=\textwidth,trim=2 2 2 2,clip]{./F_measure_visc_mms/figures/figs_res1_vs_cpu_g04.jpeg}
          \caption{Grid10.}
      \end{subfigure}
            \caption{
\label{fig:mms_res_cpu_glvl04}%
Grid level 4: Residual norm versus CPU time for the MMS test case.
} 
\end{figure}

\begin{figure}[htbp!]
    \centering  
          \begin{subfigure}[t]{0.24\textwidth}
        \includegraphics[width=\textwidth,trim=2 2 2 2,clip]{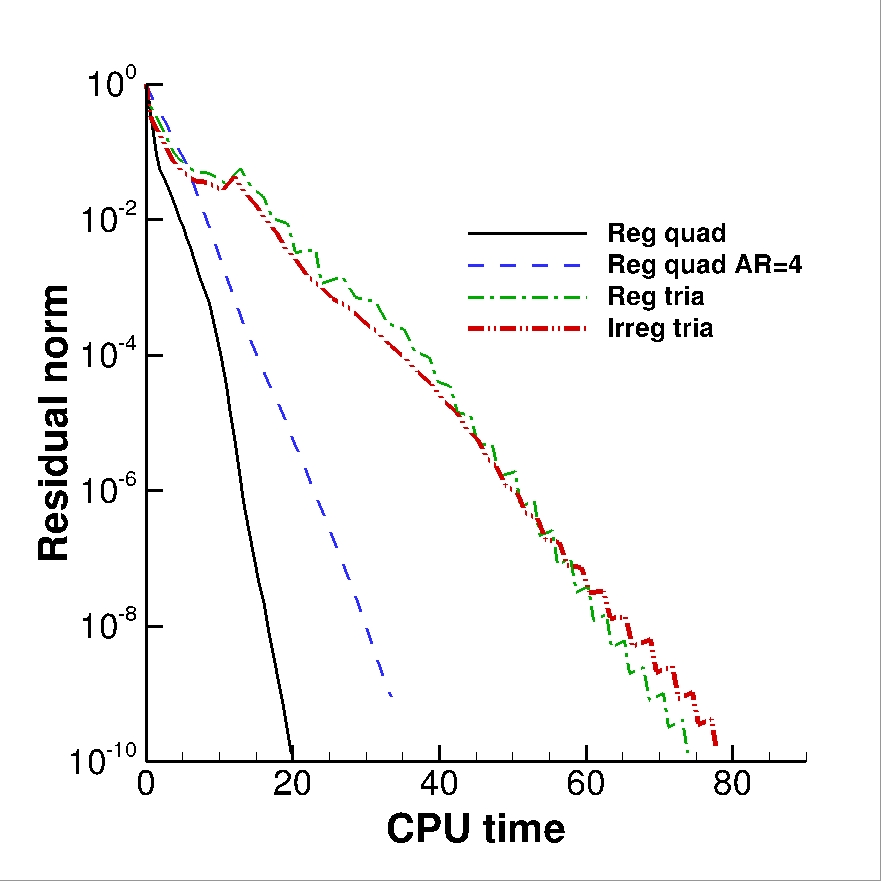}
          \caption{Grid00.}
      \end{subfigure}
          \hfill    
          \begin{subfigure}[t]{0.24\textwidth}
        \includegraphics[width=\textwidth,trim=2 2 2 2,clip]{./F_measure_visc_mms/figures/figs_res1_vs_cpu_g05.jpeg}
          \caption{Grid01.}
      \end{subfigure}
          \begin{subfigure}[t]{0.24\textwidth}
        \includegraphics[width=\textwidth,trim=2 2 2 2,clip]{./F_measure_visc_mms/figures/figs_res1_vs_cpu_g05.jpeg}
          \caption{Grid02.}
      \end{subfigure}
          \hfill    
          \begin{subfigure}[t]{0.24\textwidth}
        \includegraphics[width=\textwidth,trim=2 2 2 2,clip]{./F_measure_visc_mms/figures/figs_res1_vs_cpu_g05.jpeg}
          \caption{Grid10.}
      \end{subfigure}
            \caption{
\label{fig:mms_res_cpu_glvl05}%
Grid level 5: Residual norm versus CPU time for the MMS test case.
} 
\end{figure}

Iterative convergence results are shown for all grid levels in Figures \ref{fig:mms_res_cpu_glvl01}-\ref{fig:mms_res_cpu_glvl05}.
In all cases, as expected from the G-measure, the solver converges fast in the order for Grid00, Grid01, Grid 03, and then Grid10.
Note that the CPU time is normalized by the average CPU time spent on one residual 
computation during the iteration, which is necessary since the grid size varies.

%%%%%%%%%%%%%%%%%%%%%%%%%%%%%%%%%
\subsection{Viscous flow over a Joukowsky airfoil at $M_\infty = 0.5$ and $Re_\infty = 1000$}
%%%%%%%%%%%%%%%%%%%%%%%%%%%%%%%%%
%%%%%%%%%%%%%%%%%%%%%%%%%%%%%%%%%

\begin{figure}[htbp!]
    \centering
          \begin{subfigure}[t]{0.39\textwidth}
        \includegraphics[width=\textwidth,trim=2 2 2 2,clip]{./F_measure_inv_airfoil_tria/case_lsq_grid00/figs_grid.jpeg}
                  \caption{Grid00.}
          \label{fig:airfoil_grid00}
      \end{subfigure}
    \hfill    
          \begin{subfigure}[t]{0.39\textwidth}
        \includegraphics[width=\textwidth,trim=2 2 2 2,clip]{./F_measure_inv_airfoil_tria/case_lsq_grid01/figs_grid.jpeg}
                  \caption{Grid01.}
          \label{fig:airfoil_grid01}
      \end{subfigure}
    \hfill    
          \begin{subfigure}[t]{0.39\textwidth}
        \includegraphics[width=\textwidth,trim=2 2 2 2,clip]{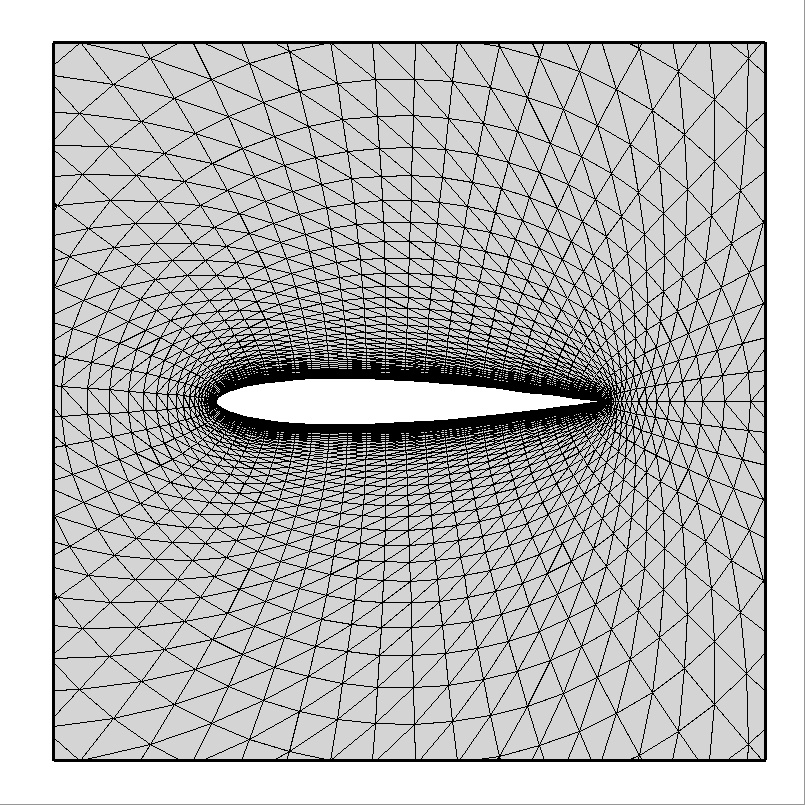}
                  \caption{Grid02.}
          \label{fig:airfoil_grid02}
      \end{subfigure}
    \hfill    
          \begin{subfigure}[t]{0.39\textwidth}
        \includegraphics[width=\textwidth,trim=2 2 2 2,clip]{./F_measure_inv_airfoil_tria/case_lsq_grid10/figs_grid.jpeg}
                  \caption{Grid10.}
          \label{fig:airfoil_grid10}
      \end{subfigure}
            \caption{
\label{fig:airfoil_grids}%
Grids used for the laminar Joukowsky airfoil test case.
} 
\end{figure}

\begin{figure}[htbp!]
    \centering
          \begin{subfigure}[t]{0.4\textwidth}
        \includegraphics[width=\textwidth,trim=2 2 2 2,clip]{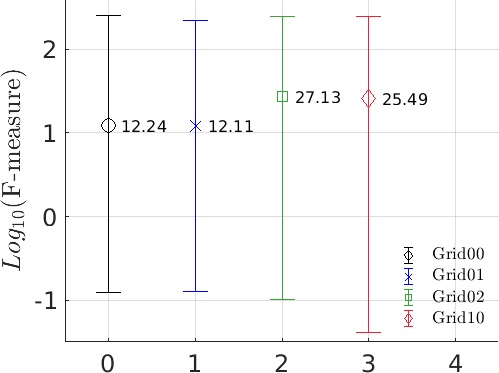}
          \caption{F-measure.}
          \label{fig:airfoil_f}
      \end{subfigure}
          \hfill    
          \begin{subfigure}[t]{0.4\textwidth}
        \includegraphics[width=\textwidth,trim=2 2 2 2,clip]{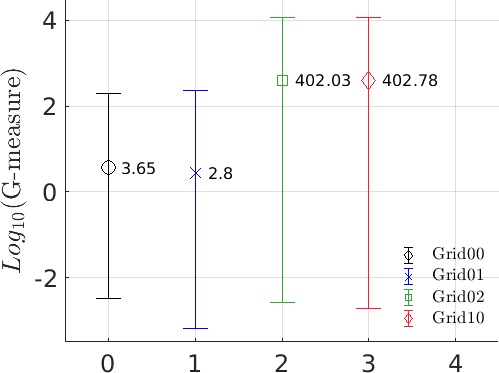}
                  \caption{G-measure.}
          \label{fig:airfoil_g}
      \end{subfigure}
                  \caption{
\label{fig:airfoil_f_and g}%
F- and G-measures for the Joukowsky airfoil grids.
} 
      \hfill    
\end{figure}

\begin{figure}[htbp!]
    \centering  
          \begin{subfigure}[t]{0.48\textwidth}
        \includegraphics[width=\textwidth,trim=2 2 2 2,clip]{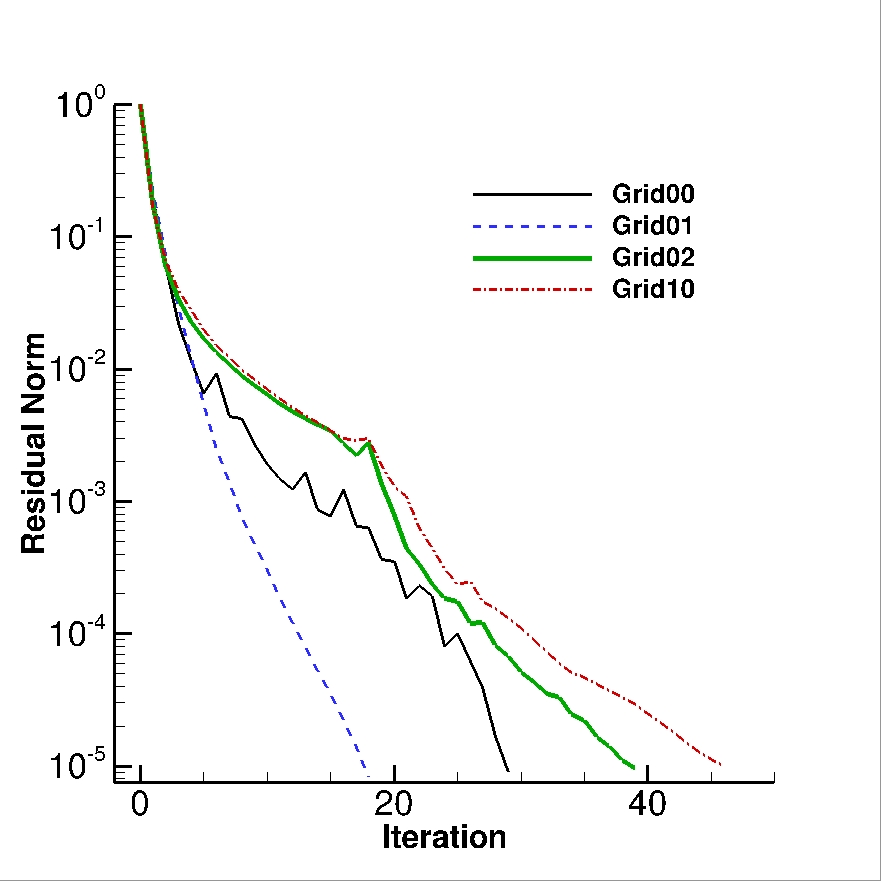}
          \caption{Residual versus iteration.}
          \label{fig:airfoil_itr}
      \end{subfigure}
          \hfill    
          \begin{subfigure}[t]{0.48\textwidth}
        \includegraphics[width=\textwidth,trim=2 2 2 2,clip]{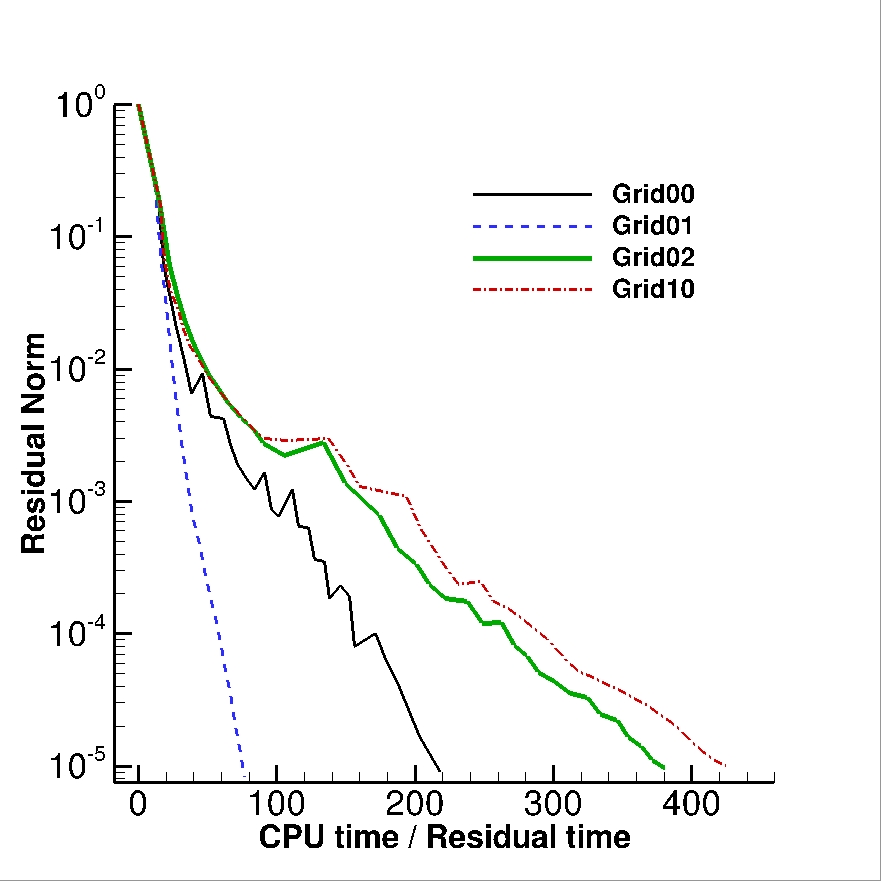}
          \caption{Residual versus normalized CPU time.}
          \label{fig:airfoil_cpu}
      \end{subfigure}
            \caption{
\label{fig:airfoil_conv}%
Convergence results for the laminar Joukowsky airfoil test case.
} 
\end{figure}

We consider solving a viscous flow over a Joukowsky airfoil at $M_\infty = 0.5$ and $Re_\infty = 1000$ with $T_\infty = 300 [K]$
and the Prandtl number 0.72, on the four different grids shown in Figure \ref{fig:airfoil_grids}. 
The F- and G-measures for these grids are computed and compared in Figure \ref{fig:airfoil_f_and g}.
As can be seen, we have for the F-measure:
\begin{eqnarray}
\mbox{Grid01} < \mbox{Grid00} < \mbox{Grid10} < \mbox{Grid02},
\end{eqnarray}
and for the G-measure:
\begin{eqnarray}
\mbox{Grid01} < \mbox{Grid00} < \mbox{Grid02} < \mbox{Grid10}.
\end{eqnarray}
The convergence results shown in Figure \ref{fig:airfoil_conv} indicate that 
the G-measure gives a more accurate indication for convergence behavior.
See Figure \ref{fig:airfoil_cpu}. 
Note again that the CPU time is normalized by the average CPU time spent on one residual 
computation during the iteration, which is necessary since the grid size varies.

%%%%%%%%%%%%%%%%%%%%%%%%%%%%%%%%%
%%%%%%%%%%%%%%%%%%%%%%%%%%%%%%%%%
\subsection{Viscous flow over a flat plate at $M_\infty = 0.15$ and $Re_\infty = 10^6$}
%%%%%%%%%%%%%%%%%%%%%%%%%%%%%%%%%
%%%%%%%%%%%%%%%%%%%%%%%%%%%%%%%%%

\begin{figure}[htbp!]
    \centering
        \includegraphics[width=0.45\textwidth,trim=2 2 2 2,clip]{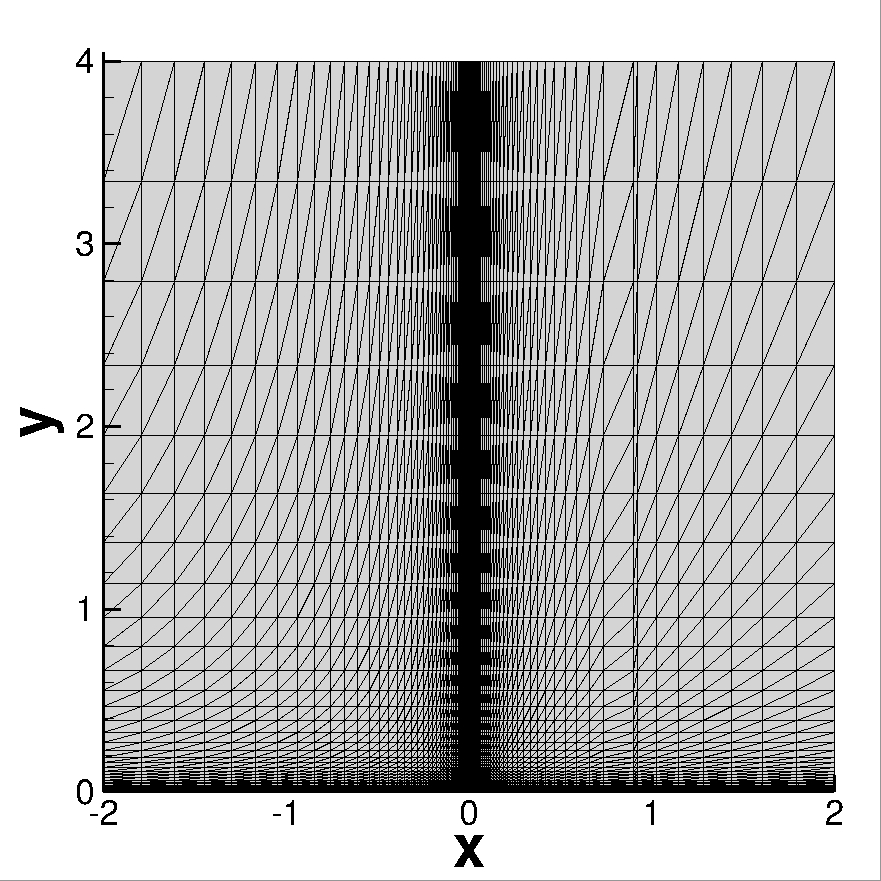}
            \caption{
\label{fig:fp_grid00_entire}%
Grid00 for the flat plate test case.
} 
\end{figure}

\begin{figure}[htbp!]
    \centering
          \begin{subfigure}[t]{0.24\textwidth}
  \centering
        \includegraphics[width=\textwidth,trim=2 2 2 2,clip]{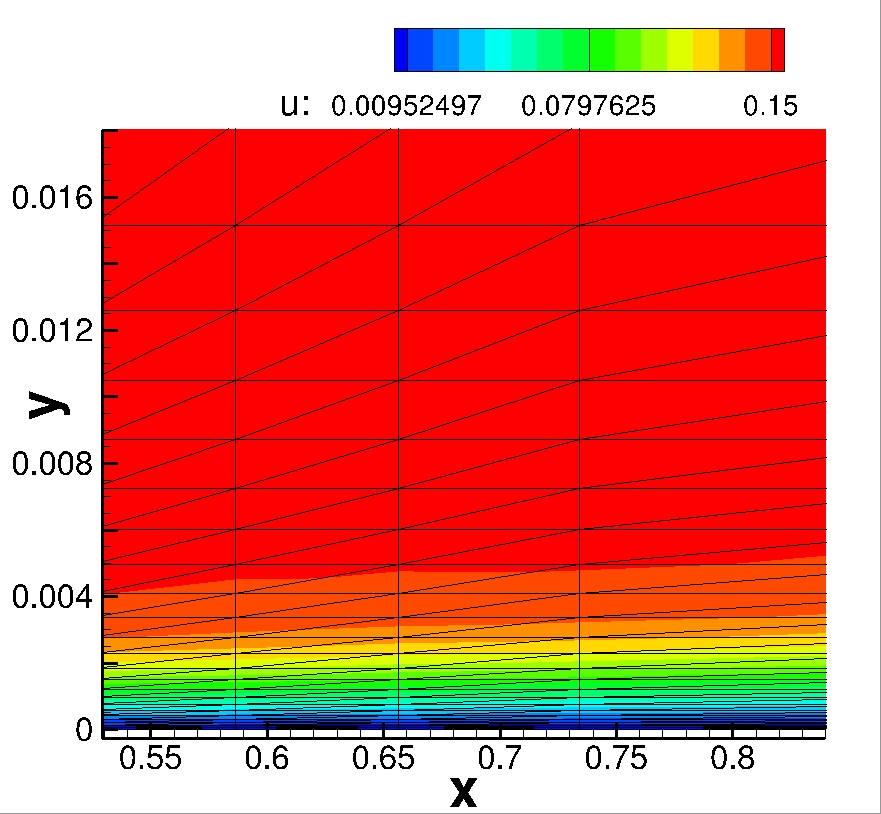}
                  \caption{Grid00.}
          \label{fig:fp_grid00}
      \end{subfigure}
      \hfill    
          \begin{subfigure}[t]{0.24\textwidth}
  \centering
        \includegraphics[width=\textwidth,trim=2 2 2 2,clip]{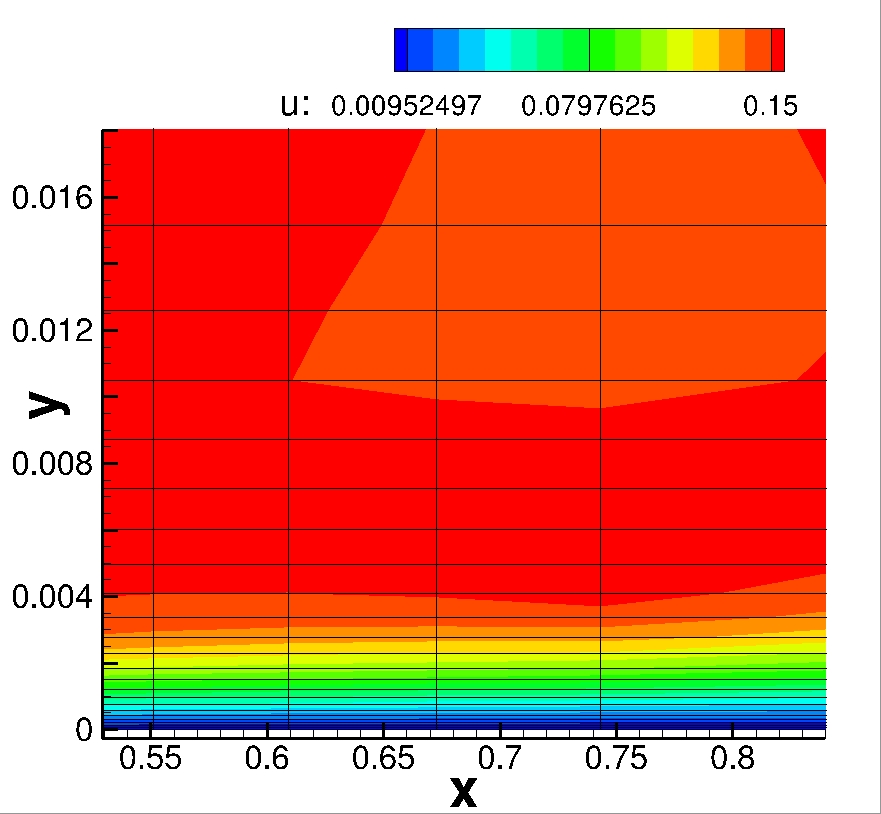}
                  \caption{Grid01.}
          \label{fig:fp_grid01}
      \end{subfigure}
          \begin{subfigure}[t]{0.24\textwidth}
  \centering
        \includegraphics[width=\textwidth,trim=2 2 2 2,clip]{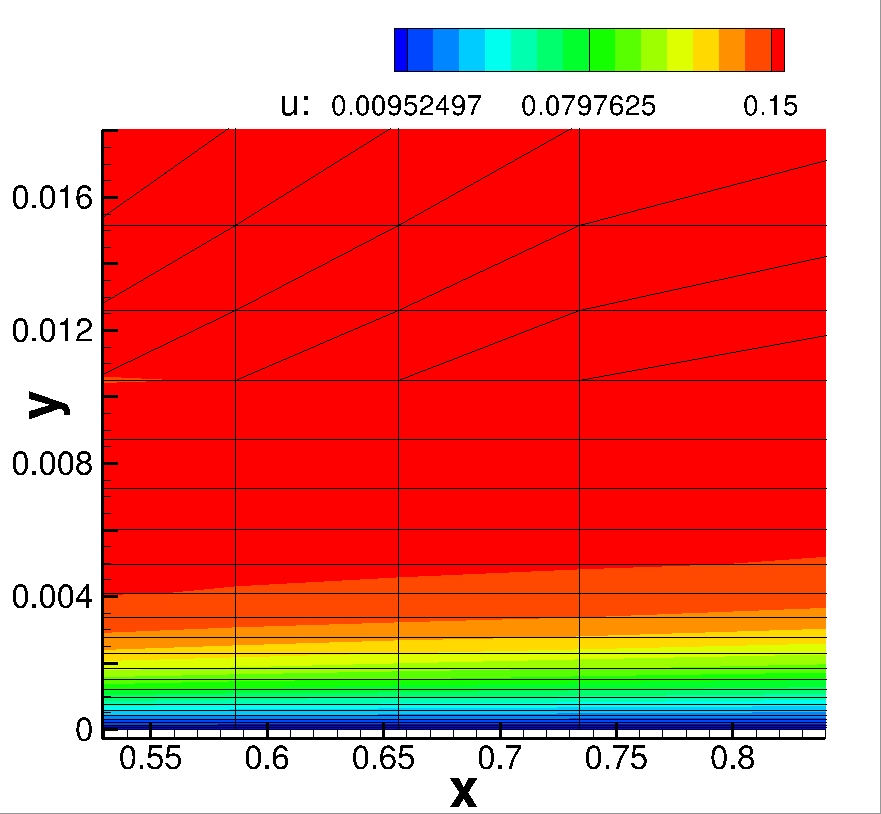}
                  \caption{Grid02.}
          \label{fig:fp_grid00}
      \end{subfigure}
      \hfill    
          \begin{subfigure}[t]{0.24\textwidth}
  \centering
        \includegraphics[width=\textwidth,trim=2 2 2 2,clip]{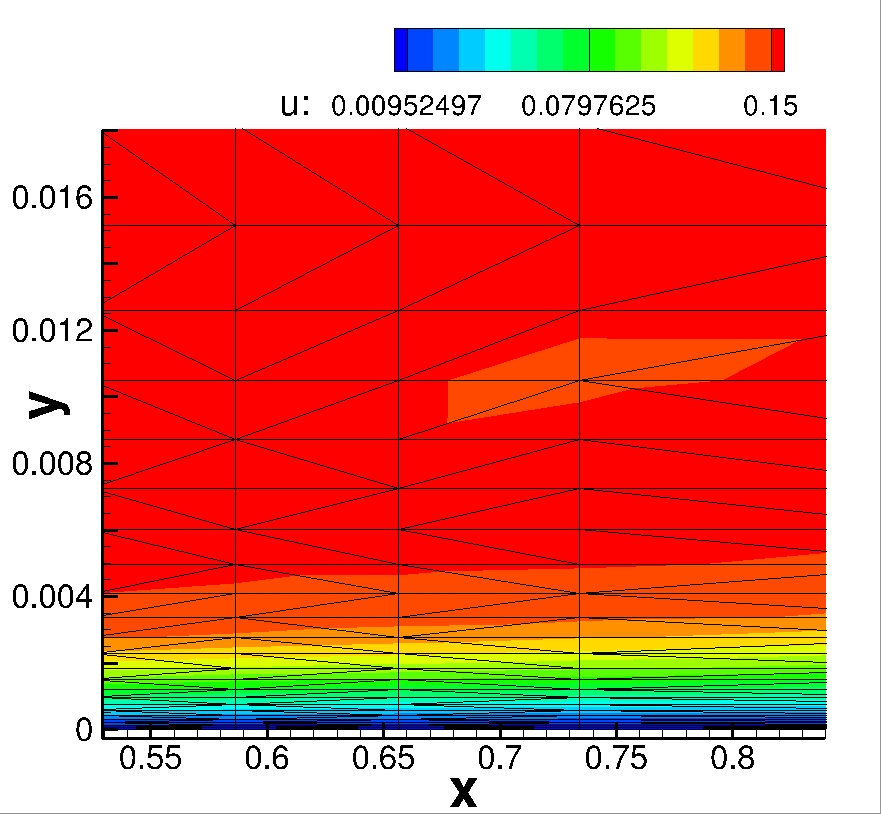}
                  \caption{Grid10.}
          \label{fig:fp_grid01}
      \end{subfigure}
                  \caption{
\label{fig:cylinder_grid}%
Grids used for the flat plate test case.
}
\end{figure}

\begin{figure}[htbp!]
    \centering
          \begin{subfigure}[t]{0.4\textwidth}
        \includegraphics[width=\textwidth,trim=2 2 2 2,clip]{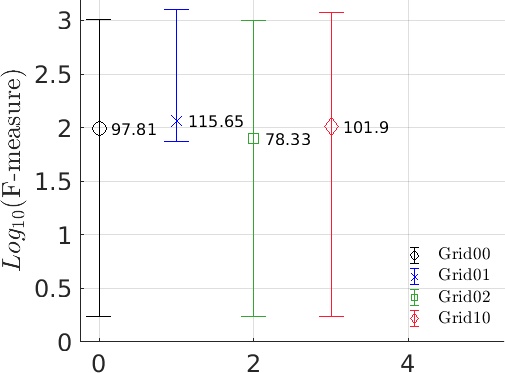}
          \caption{F-measure.}
          \label{fig:airfoil_irreg_centroids}
      \end{subfigure}
          \hfill    
          \begin{subfigure}[t]{0.4\textwidth}
        \includegraphics[width=\textwidth,trim=2 2 2 2,clip]{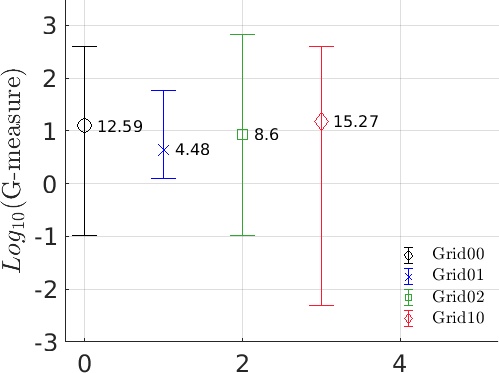}
          \caption{G-measure.}
          \label{fig:airfoil_irreg_skewness}
      \end{subfigure}
                  \caption{
\label{fig:fp_f_and_g}%
F- and G-measures for the flat plate grids.
} 
      \hfill    
\end{figure}

\begin{figure}[htbp!]
    \centering
          \begin{subfigure}[t]{0.48\textwidth}
        \includegraphics[width=\textwidth,trim=2 2 2 2,clip]{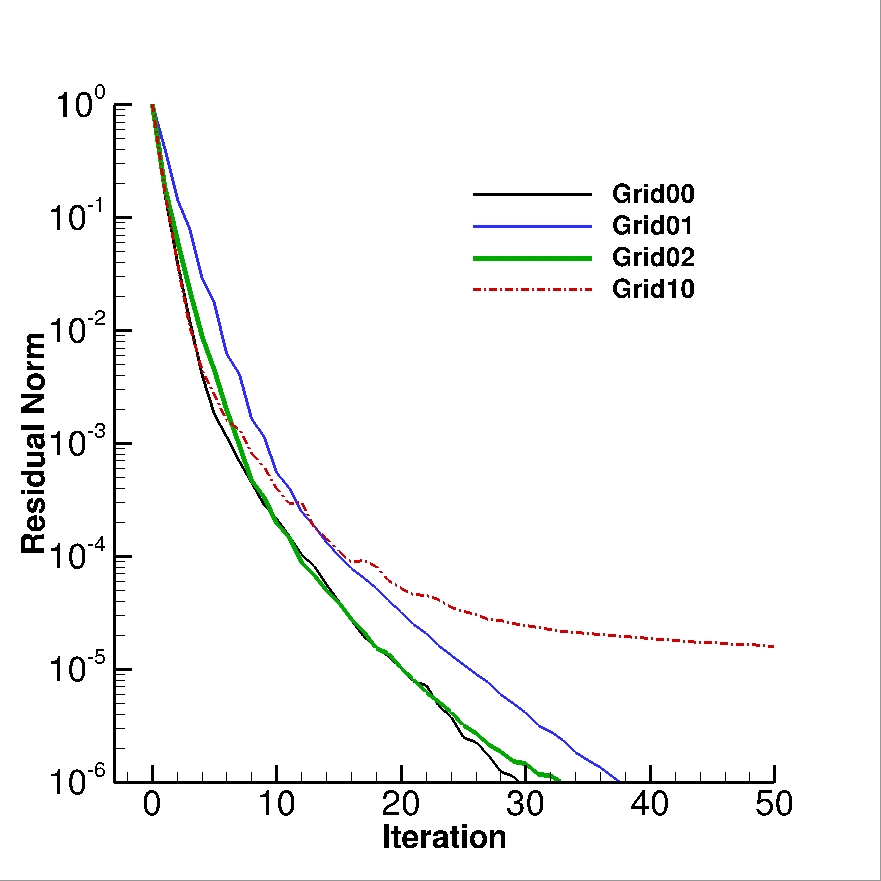}
          \caption{Residual versus iteration.}
          \label{fig:fp_itr}
      \end{subfigure}
      \hfill    
          \begin{subfigure}[t]{0.48\textwidth}
        \includegraphics[width=\textwidth,trim=2 2 2 2,clip]{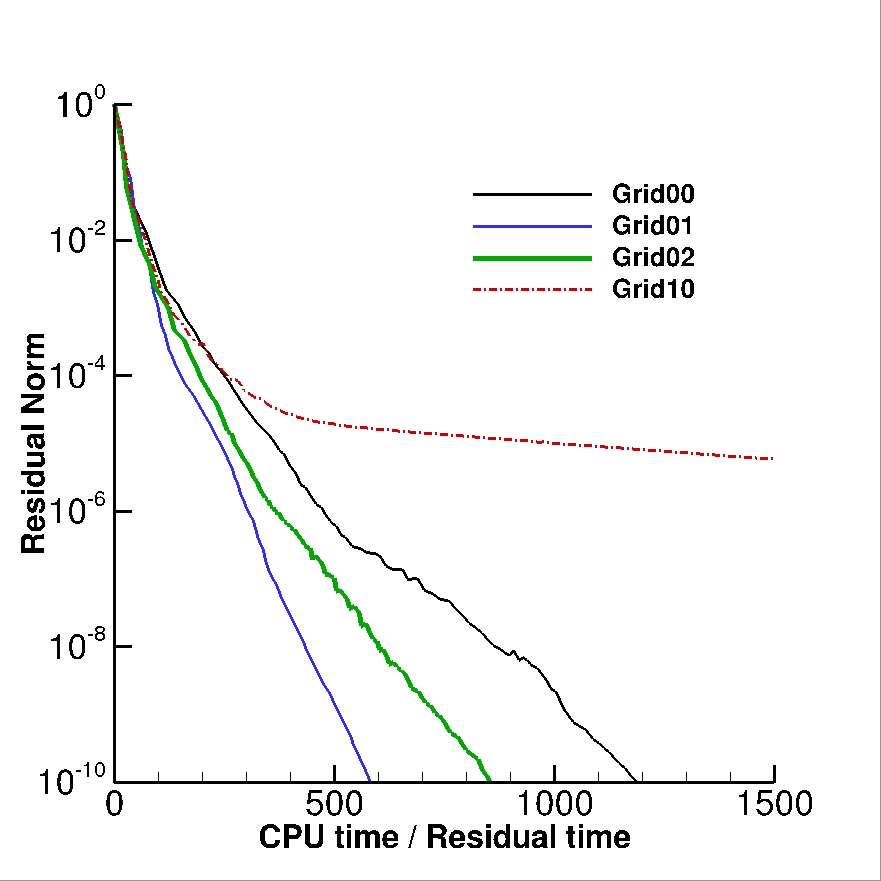}
          \caption{Residual versus normalized CPU time.}
          \label{fig:fp_cpu}
      \end{subfigure}
            \caption{
\label{fig:fp_conv}%
Convergence results for the flat plate test case.
} 
\end{figure}

Another test was performed for a viscous flow over a flat plate at the free stream Mach number 
$M_\infty = 0.15$ and $Re_\infty = 10^6$, where $Re_\infty$ is the free stream Reynolds number 
per grid unit. The domain is a rectangular definedby $[x,y] = [-2,2] \times [0,4]$ with the flat plate in
$x = [0,2]$ at $y=0$. See Figure \ref{fig:fp_grid00_entire}.
 Four different grids are compared. See Figure \ref{fig:cylinder_grid}. 
The F- and G-measures for these grids are computed and compared in Figure \ref{fig:fp_f_and_g}.
As can be seen, we have for the F-measure:
\begin{eqnarray}
\mbox{Grid00} < \mbox{Grid02} < \mbox{Grid10} < \mbox{Grid01},
\end{eqnarray}
and for the G-measure:
\begin{eqnarray}
\mbox{Grid01} < \mbox{Grid02} < \mbox{Grid00} < \mbox{Grid10}.
\end{eqnarray}
Again, the G-measure is consistent with the actual convergence results.
See Figure \ref{fig:fp_conv}.

%%%%%%%%%%%%%%%%%%%%%%%%%%%%%%%%%%%%%%%%%%%%%%%%%%%%%%%%%%%
\section{Concluding Remarks}
%%%%%%%%%%%%%%%%%%%%%%%%%%%%%%%%%%%%%%%%%%%%%%%%%%%%%%%%%%%

We have proposed two grid-quality measures, F-measure and G-measure, in relation to iterative convergence of an implicit 
finite-volume solver for unstructured grids. The F-measure is a lower bound of a LSQ gradient, that can be computed
for a given grid. The G-measure is the magnitude of a LSQ gradient of a known function whose gradient is zero a cell centroid.
These measures are defined specifically for unstructured-grid finite-volume solvers, employing LSQ gradient methods.
Numerical results indicate that the G-measure more accurately the performance of the implicit solver: it converges faster on a grid with 
a smaller G-measure. 

In future work, the G-measure should be tested for practical problems with three-dimensional unstructured grids and also 
for node-centered methods such as the edge-based discretization. 
Furthermore, it should be investigated whether a grid can be modified to minimize the G-measure.
For example, a node may be removed/added or an edge may be swapped if it reduces the G-measure locally.

%%%%%%%%%%%%%%%%%%%%%%%%%%%%%%%%%%%%%%%%
\addcontentsline{toc}{section}{Acknowledgments}
\section*{Acknowledgments}

The author gratefully acknowledges support from Software CRADLE, part of Hexagon.

%%%%%%%%%%%%%%%%%%%%%%%%%%%%%%%%%%%%%%%%%%%%%%%%%%%%%%%%%%%
\bibliography{../../bibtex_nishikawa_database}

\begin{thebibliography}{10}

\bibitem{VeluriRoyLuke:CF2012}
S.~P. Veluri, C.~J. Roy, and E.~A. Luke.
\newblock Comprehensive code verification techniques for finite volume {CFD}
  codes.
\newblock {\em Comput. Fluids}, 70:59--72, 2012.

\bibitem{KallinderisFotia:JCP2015}
Y.~Kallinderis and S.~Fotia.
\newblock A priori mesh quality metrics for three-dimensional hybrid grids.
\newblock {\em J. Comput. Phys.}, 280:465--488, 2015.

\bibitem{YouMittalWangMoin:JCP2006}
D.~You, R.~Mittal, M.~Wang, and P.~Moin.
\newblock Analysis of stability and accuracy of finite-difference schemes on a
  skewed mesh.
\newblock {\em J. Comput. Phys.}, 213:184--204, 2006.

\bibitem{DennerWachem:JCP2015}
F.~Denner and B.~{G}.~{M}. van Wachem.
\newblock {TVD} differencing on three-dimensional unstructured meshes with
  monotonicity-preserving correction of mesh skewness.
\newblock {\em J. Comput. Phys.}, 298:466--479, 2015.

\bibitem{ghoreyshi_etal:AIAA2015-0407}
M.~Ghoreyshi, K.~Bergeron, J.~Seidel, An.~{J}. Lofthouse, and R.~M. Cummings.
\newblock Grid quality and resolution effects on the aerodynamic modeling of
  ram-air parachute canopies.
\newblock In {\em Proc. of 53rd AIAA Aerospace Sciences Meeting}, {AIAA} Paper
  2015-0407, Kissimmee, Florida, January 2015.

\bibitem{nc_cc_comparison_viscous:AIAAJ2010}
B.~Diskin, J.~L. Thomas, E.~J. Nielsen, H.~Nishikawa, and J.~A. White.
\newblock Comparison of node-centered and cell-centered unstructured
  finite-volume discretizations: {V}iscous fluxes.
\newblock {\em {AIAA} J.}, 48(7):1326--1338, July 2010.

\bibitem{Dannenhoffer:AIAA2012-0610}
J.~Dannenhoffer.
\newblock Correlation of grid quality metrics and solution accuracy for
  supersonic flows.
\newblock In {\em Proc. of 50th AIAA Aerospace Sciences Meeting}, {AIAA} Paper
  2012-0610, Nashville, Tennessee, 2012.

\bibitem{pointwise_article_link}
{P}ointwise John~Rhoads.
\newblock How do you define a good grid?
\newblock
  \url{https://www.pointwise.com/theconnector/2014-November/How-Do-You-Define-Good-Grid.html},
  visited on 01/23/2019.

\bibitem{diskin_thomas:AIAAJ2011}
B.~Diskin and J.~L. Thomas.
\newblock Comparison of node-centered and cell-centered unstructured
  finite-volume discretizations: {I}nviscid fluxes.
\newblock {\em {AIAA} J.}, 49(4):836--854, 2011.

\bibitem{nishikawa_nakashima_watanabe:jcp2017}
H.~Nishikawa, Y.~Nakashima, and N.~Watanabe.
\newblock Effects of high-frequency damping on iterative convergence of
  implicit viscous solver.
\newblock {\em J. Comput. Phys.}, 348:66--81, 2017.

\bibitem{nishikawa_stencil:JCP2019}
H.~Nishikawa.
\newblock Efficient gradient stencils for robust implicit finite-volume solver
  convergence on distorted grids.
\newblock {\em J. Comput. Phys.}, 386:486--501, 2019.

\bibitem{Kleb_etal_aiaa2019-2948}
W.~{L}. Kleb, M.~{A}. Park, W.~{A}. Wood, K.~{L}. Bibb, K.~{B}. Thompson, and
  R.~{J}. Gomez.
\newblock Sketch-to-solution: An exploration of viscous {CFD} with automatic
  grids.
\newblock In {\em 24th {AIAA} Computational Fluid Dynamics Conference}, {AIAA}
  Paper 2019-2948, Dallas, TX, 2019.

\bibitem{WhiteNishikawaBaurle_aiaa2019-0127}
J.~{A}. White, H.~Nishikawa, and R.~Baurle.
\newblock Weighted least-squares cell-average gradient construction methods for
  the {VULCAN-CFD} second-order accurate unstructured grid cell-centered
  finite-volume solver.
\newblock In {\em AIAA Scitech 2019 Forum}, {AIAA} Paper 2019-0127, San
  {D}iego, CA, 2019.

\bibitem{nishikawa_centroid:JCP2020}
H.~Nishikawa.
\newblock A face-area-weighted centroid formula for finite-volume method that
  improves skewness and convergence on triangular grids.
\newblock {\em J. Comput. Phys.}, 401:109001, 2020.

\end{thebibliography}
\bibliographystyle{unsrt}
%%%%%%%%%%%%%%%%%%%%%%%%%%%%%%%%%%%%%%%%%%%%%%%%%%%%%%%%%%%

%%%%%%%%%%%%%%%%%%%%%%%%%%%%%%%%%%%%%%%%
%%%%%%%%%%%%%%%%%%%%%%%%%%%%%%%%%%%%%%%%%

%%%%%%%%%%%%%%%%%%%%%%%%%%%%%%%%%%%%%%%%
%%%%%%%%%%%%%%%%%%%%%%%%%%%%%%%%%%%%%%%%%

%%%%%%%%%%%%%%%%%%%%%%%%%%%%%%%%%%%%%%%%
%%%%%%%%%%%%%%%%%%%%%%%%%%%%%%%%%%%%%%%%%
%%%%%%%%%%%%%%%%%%%%%%%%%%%%%%%%%%%%%%%%
%%%%%%%%%%%%%%%%%%%%%%%%%%%%%%%%%%%%%%%%%
%%%%%%%%%%%%%%%%%%%%%%%%%%%%%%%%%%%%%%%%
%%%%%%%%%%%%%%%%%%%%%%%%%%%%%%%%%%%%%%%%%
%%%%%%%%%%%%%%%%%%%%%%%%%%%%%%%%%%%%%%%%
%%%%%%%%%%%%%%%%%%%%%%%%%%%%%%%%%%%%%%%%%
%%%%%%%%%%%%%%%%%%%%%%%%%%%%%%%%%%%%%%%%
%%%%%%%%%%%%%%%%%%%%%%%%%%%%%%%%%%%%%%%%%

%%%%%%%%%%%%%%%%%%%%%%%%%%%%%%%%%%%%%%%%
%%%%%%%%%%%%%%%%%%%%%%%%%%%%%%%%%%%%%%%%%

%%%%%%%%%%%%%%%%%%%%%%%%%%%%%%%%%%%%%%%%
%%%%%%%%%%%%%%%%%%%%%%%%%%%%%%%%%%%%%%%%%
%%%%%%%%%%%%%%%%%%%%%%%%%%%%%%%%%%%%%%%%
%%%%%%%%%%%%%%%%%%%%%%%%%%%%%%%%%%%%%%%%%

%%%%%%%%%%%%%%%%%%%%%%%%%%%%%%%%%%%%%%%%
%%%%%%%%%%%%%%%%%%%%%%%%%%%%%%%%%%%%%%%%%

%%%%%%%%%%%%%%%%%%%%%%%%%%%%%%%%%%%%%%%%
%%%%%%%%%%%%%%%%%%%%%%%%%%%%%%%%%%%%%%%%%

%%%%%%%%%%%%%%%%%%%%%%%%%%%%%%%%%%%%%%%%
%%%%%%%%%%%%%%%%%%%%%%%%%%%%%%%%%%%%%%%%%

%%%%%%%%%%%%%%%%%%%%%%%%%%%%%%%%%%%%%%%%%
%%%%%%%%%%%%%%%%%%%%%%%%%%%%%%%%%%%%%%%%
%%%%%%%%%%%%%%%%%%%%%%%%%%%%%%%%%%%%%%%%%

\end{document}